\documentclass[12pt]{amsart}
\usepackage{amscd}

\newcommand{\C}{{\mathbb{C}}}        
\newcommand{\hyp}{\mathcal{H}}

\newcommand{\G}{{\mathcal{G}}}       

\newcommand{\oddots}{{\mathinner{\mkern1mu\raise1pt\vbox{\kern7pt\hbox{.}}\mkern2mu\raise4pt\hbox{.}\mkern2mu\raise7pt\hbox{.}\mkern1mu}}}

\newcommand{\As}{{(A,\sigma)}}

\newcommand{\Fx}{{F^{\ast}}}

\newcommand{\supast}[1]{{{#1}^\ast}}
\newcommand{\xsq}[1]{{{#1}^{\ast2}}}
\newcommand{\sq}[1]{{\supast{#1}/\xsq{#1}}}

\newcommand{\qform}[1]{{\langle{#1}\rangle}}                   

\newcommand{\basemu}{\mathbf{\mu}}
\newcommand{\mmu}[1]{\basemu_#1}     

\newcommand{\Skew}{\mathrm{Skew}}

\newcommand{\Int}{{\mathrm{Int}\,}}

\newcommand{\Gal}{{\mathrm{{\mathcal{G}}a\ell}\,}}
\newcommand{\im}{{\mathrm{im}\,}}

\DeclareMathOperator{\Id}{Id}

\newcommand{\End}[1]{{\mathrm{End}_{#1}}}
\newcommand{\EndF}{\End{F}}

\newcommand{\Aut}[1]{{{\mathrm {Aut}}_{#1}}}
\newcommand{\aut}[1]{\Aut{#1}}

\newcommand{\Hom}{{\mathrm{Hom}}}

\newcommand{\osum}{\oplus}

\newcommand{\longto}{\longrightarrow}
\newcommand{\iso}{\xrightarrow{\sim}}
\newcommand{\injects}{\hookrightarrow}

\newcommand{\Gbar}{\overline{G}}

\newcommand{\wdet}{\mathop{\mathrm{wdet}}}

\newcommand{\Rsq}{\sq{R}}

\newcommand{\Ax}{\supast{A}}

\newcommand{\s}{\sigma}

\newcommand{\bbrown}{(B,{-})}
\newcommand{\bgift}{(A, \s, \pi)}

\newcommand{\tbtmat}[4]{\left( \begin{array}{cc} #1&#2 \\ #3&#4 \end{array} \right) }
\newcommand{\stbtmat}[4]{\left( \begin{smallmatrix} #1&#2 \\ #3&#4 \end{smallmatrix} \right) }
\newcommand{\thbthmat}[9]{\left( \begin{array}{ccc} #1&#2&#3 \\ #4&#5&#6 \\ #7&#8&#9 \end{array} \right) }
\newcommand{\sthbthmat}[9]{\left( \begin{smallmatrix} #1&#2&#3 \\ #4&#5&#6 \\ #7&#8&#9 \end{smallmatrix} \right) }

\newcommand{\basmat}{\tbtmat{\alpha}{j}{j'}{\beta}}
\newcommand{\EndK}{\mathrm{End}_{K}\,}
\newcommand{\HomF}{\mathrm{Hom}_{F}\,}

\newcommand{\tr}{\mathrm{tr}}

\newcommand{\st}{\s_2}
\newcommand{\Sand}{\mathrm{Sand}}

\newcommand{\pihat}{\widehat{\pi}}
\newcommand{\shat}{\widehat{\s}}
\newcommand{\idhat}{\widehat{\Id}}

\newcommand{\Bq}{\B^q}
\newcommand{\Md}{\M^d}
\newcommand{\Ms}{\M_s}
\newcommand{\EndM}{\mathrm{End}(\M)}
\newcommand{\EndB}{\mathrm{End}(\B)}

\newcommand{\Der}{\mathrm{Der}\,}

\newcommand{\TrdA}{\mathop{\mathrm{Trd}_A}}

\newcommand{\gr}[2]{^{#1}\hskip-.2em #2}

\newcommand{\Inv}{\mathrm{Inv}\,}
\newcommand{\Iso}{\mathrm{Iso}\,}

\newcommand{\autp}{\Aut{}^{+}\,}
\renewcommand{\aut}{\Aut{}\,}
\renewcommand{\Sp}{\mathrm{Sp}\,}
\newcommand{\PSp}{\aut}
\newcommand{\GmF}{\mathbb{G}_{m,F}}

\newcommand{\B}{\mathcal{B}}

\newcommand{\h}{\mathfrak{H}}
\newcommand{\gift}{\mathfrak{G}}
\renewcommand{\C}{\mathcal{C}}

\renewcommand{\G}{\gift}
\newcommand{\M}{\mathfrak{M}}
\newcommand{\D}{\Delta}

\newtheorem{thm}{Theorem}[section]        
\newtheorem{lem}[thm]{Lemma}

\newtheorem{prop}[thm]{Proposition}

\newtheorem{const}[thm]{Construction}

\theoremstyle{definition}

\newtheorem{defn}[thm]{Definition}
\newtheorem{eg}[thm]{Example}

\theoremstyle{remark}

\newtheorem{rmk}[thm]{Remark}

\newenvironment{pf}{\noindent{\sl Proof:}}{\hfill$\qed$\medskip}

\numberwithin{equation}{section}

\makeatletter
\newenvironment{neqn}%
{\setcounter{equation}{\value{thm}}\begin{eqnarray}}%
{\end{eqnarray}{\stepcounter{thm}}\global\@ignoretrue}
\makeatother

\makeatletter
{\setcounter{equation}{\value{thm}}}%
{{\stepcounter{thm}}\global\@ignoretrue}
\makeatother

\makeatletter
{\setcounter{equation}{\value{thm}}}%
{{\stepcounter{thm}}\global\@ignoretrue}
\makeatother

\makeatletter
{\smallskip \noindent\refstepcounter{thm}{\bf \thethm.}{\bf{ #1.}}}%
{\smallskip \global\@ignoretrue}
\makeatother

\makeatletter
\newenvironment{borel*}%
{\smallskip \noindent\refstepcounter{thm}{\bf \thethm.}}%
{\smallskip \global\@ignoretrue}
\makeatother

\begin{document}

\title{Groups of type $E_7$ over arbitrary fields}
\author{R.~Skip Garibaldi}


\begin{abstract}
Freudenthal triple systems come in two flavors, degenerate and
nondegenerate.  The best criterion for distinguishing between the two
which is available in the literature is by descent.  We provide an identity
which is satisfied only by nondegenerate triple systems.  We then use this to
define algebraic structures whose automorphism groups produce all
adjoint algebraic 
groups of type $E_7$ over an arbitrary field of characteristic $\ne 2,
3$. 

The main advantage of these new structures is that they incorporate a
previously unconsidered invariant (a symplectic involution) of these
groups in a fundamental way.
As an application, we give a construction of adjoint groups with
Tits algebras of index 2 which provides a complete description of
this involution and apply this to groups of type $E_7$
over a real-closed field.
\end{abstract}

\maketitle

A useful strategy for studying simple (affine) algebraic groups over
arbitrary fields has been to describe such a group as the group of
automorphisms of some algebraic object.  We restrict our attention to
fields of characteristic $\ne 2, 3$.  The idea is that these 
algebraic objects are easier to study, and their properties
correspond to properties of the group one is interested in.  Weil
described all groups of type $A_n$, $B_n$, $C_n$, $\gr{1}{D_n}$, and
$\gr{2}{D_n}$ 
in this manner in
\cite{Weil}.  Similar descriptions were soon found for groups of type
$F_4$ (as automorphism groups of Albert algebras) and $G_2$ (as
automorphism groups 
of Cayley algebras).  Recently, groups of type $\gr{3}{D_4}$ and $\gr{6}{D_4}$
have been described in \cite[\S 43]{KMRT} as groups of automorphisms of
trialitarian central simple algebras.  The remaining groups are
those of types $E_6$, $E_7$, and $E_8$.  

As an attempt to provide an algebraic structure associated to groups
of type $E_7$, Freudenthal introduced a new kind of algebraic
structure in \cite[\S 4]{Frd:E7.1}, which was later 
studied
axiomatically in \cite{Meyb:FT}, \cite{Brown:E7}, and
\cite{Ferr:strict}.  These objects, called Freudenthal triple
systems, come in two flavors: degenerate and nondegenerate.  The
automorphism groups of the nondegenerate ones provide all simply
connected groups of type $E_7$ with trivial Tits algebras over an
arbitrary field.  (The Tits
algebras of a group $G$ are the endomorphism 
rings of certain irreducible representations of $G$.  They were
introduced in \cite{Ti:Cl}, or see \cite[\S 27.A]{KMRT} for another
treatment.) 
In fact, more is true: they are precisely the
$G$-torsors for $G$ simply connected split of type $E_7$.    

One issue that has not been addressed adequately in the study of
triple systems is how to distinguish between the two kinds.  A triple
system 
consists of a 56-dimensional vector space endowed with a nondegenerate 
skew-symmetric bilinear form and a quartic form (see \ref{ftsdef} for
a complete definition), and we say that the 
triple system  is  
nondegenerate precisely when this quartic form is irreducible when we
extend scalars to a separable closure of the base field.
There seems to be essentially no 
criterion available in the literature to distinguish
between the two types other then checking the definition.  In Section
\ref{pisec}, we show that one of 
the identities which nondegenerate triple systems are known to satisfy 
is not satisfied by degenerate ones, thus providing a means for
differentiating between the two types which doesn't require
enlarging the base field.

In Section \ref{giftsec}, we define algebraic structures whose groups of
automorphisms produce all groups of type $E_7$ up to isogeny.  Thanks
to the preceding result distinguishing between degenerate and
nondegenerate triple systems, no Galois descent is needed for this
definition.  We call these structures 
gifts (short for {\em g}eneral{\em i}zed {\em F}reudenthal {\em
t}riple system{\em s}).  They are triples $\bgift$ 
such that $A$ is a central simple $F$-algebra of degree 56, $\s$ is a
symplectic involution on $A$, and $\pi \!: A \longto A$ is an
$F$-linear map satisfying certain axioms (see \ref{giftdef} for a full
definition).  We also show an equivalence of
categories between the category of gifts over an arbitrary field $F$
and the category of adjoint (equivalently, simply connected) groups of 
type $E_7$ over $F$.  A description of the flag (a.k.a.~homogeneous
projective) varieties of an arbitrary group of type $E_7$ is then 
easily derived in Section \ref{flagsec}.

The main strength of these gifts is that they include the involution
$\s$ in an intrinsic way.  Specifically, as mentioned above every
adjoint group $G$ of type $E_7$ is isomorphic to $\aut\bgift$ for some
gift $\bgift$.  The inclusion
of the split adjoint group of type $E_7$ in the split adjoint group of
type $C_{28}$ gives via Galois cohomology that the pair $\As$ is an
invariant of the group $G$.
We give a construction (in \ref{quatconst}) of gifts with   
algebra component $A$ of index 2 (which is equivalent to constructing
groups of type $E_7$ with Tits algebra of index 2 and corresponds to
Tits' construction of analogous Lie algebras in \cite{Ti:const}, or
see \cite[\S 10]{Jac:ex}) and we  
describe the  
involution $\s$ explicitly in this case.  In the final section, we use
this construction to prove some facts about simple groups of type
$E_7$ over a real-closed field.

\subsection*{Notations and conventions}
All fields that we consider will have characteristic $\ne 2, 3$.

For $g$ an element in a group $G$, we write $\Int(g)$ for the
automorphism of $G$ given by $x \mapsto g x g^{-1}$.

For $X$ a variety over a field $F$ and $K$ any field extension of $F$,
we write $X(K)$ for the $K$-points of $X$.

When we say that an affine algebraic group (scheme) $G$ is {\em
  simple}, we mean that it is absolutely almost simple in the usual
sense (i.e., $G(F_s)$ has a finite center and no noncentral normal
subgroups for $F_s$ a separable closure of the ground field).

We write $\GmF$ for the algebraic group whose $F$-points are $\Fx$ and
$\mmu{n}$ for its subgroup group of $n$th roots of unity.

We will also follow the usual conventions for Galois cohomology and write
$H^i(F, G) := H^i(\Gal(F_s/F), G(F_s))$ for $G$ any algebraic group
over $F$, and similarly for the cocycles $Z^1(F, G)$.  For more
information about Galois cohomology, see \cite{SeLF} and \cite{SeCG}.

We follow the notation in \cite{Lam} for quadratic forms.

For $I$ a right ideal in a central simple $F$-algebra $A$, we
define the {\em rank} of $I$ to be $(\dim_F I)/ \deg A$.  Thus when
$A$ is split, so that we may write $A \cong \EndF(V)$ for some
$F$-vector space $V$, $I = \HomF(V,U)$ for some subspace $U$ of $V$
and the rank of $I$ is precisely the dimension of $U$.

\section{Background on triple systems} \label{backsec}

\begin{defn} (See, for example, \cite[p.~314]{Ferr:strict} or \cite[3.1]{G:struct}) \label{ftsdef}
A {\em (simple) Freudenthal triple system} is a 3-tuple $(V, b,
t)$ such that $V$ is a 56-dimensional vector space, $b$ is a
nondegenerate skew-symmetric bilinear form on $V$, and $t$ is a
trilinear product $t \!: V \times V \times V \longto V$.

We define a 4-linear form $q(x,y,z,w) := b(x, t(y, z, w))$ for $x$, $y$, $z$, $w \in V$, and we 
require that
\begin{description}
\item[FTS1] $q$ is symmetric,
\item[FTS2] $q$ is not identically zero, and
\item[FTS3] $t(t(x,x,x),x,y) = b(y,x) t(x,x,x) + q(y, x,x,x) x$ for
all $x, y \in V$.
\end{description}
We say that such a triple system is {\em nondegenerate} if the quartic
form $v \mapsto q(v,v,v,v)$ on $V$ is absolutely irreducible (i.e.,
irreducible over a separable closure of the base field) and {\em
degenerate} otherwise.
\end{defn}

Note that since $b$ is nondegenerate, FTS1 implies that $t$ is symmetric.

One can linearize FTS3 a little bit to get an equivalent axiom that
will be of use later.  Specifically, replacing $x$ with $x + \lambda
z$, expanding using linearity, and taking the coefficient of
$\lambda^2$, one gets the equivalent formula
\[
\begin{array}{lrcl}
\mathbf{FTS3'} & t(t(x, x, z), z, y) + t(t(x, z, z), x, y) &=& z q(x, x,
z, y) + x q(x, z, z, y) \\
& && +\,b(y, z) t(x,x,z) + b(y,x) t(x, z, z).
\end{array}
\]

\begin{eg} \label{degeneg} (Cf.~\cite[p.~94]{Brown:E7},
  \cite[p.~172]{Meyb:FT}) 
Let $W$ be a 27-dimensional $F$-vector space endowed with a
non-degenerate skew-symmetric bilinear form $s$ and set
\begin{neqn} \label{matdef}
V := \tbtmat{F}{W}{W}{F}.
\end{neqn}
For 
\begin{neqn} \label{xydef}
x := \tbtmat{\alpha}{j}{j'}{\beta} \hbox{\ and\ } y :=
\tbtmat{\gamma}{k}{k'}{\delta}
\end{neqn}
set
\[
b(x, y) := \alpha \delta - \beta \gamma + s(j, k') + s(j', k).
\]
We define the {\em determinant map} $\det \!: V \longto F$ by
\[
\det(x) := \alpha \beta - s(j, j')
\]
and set
\[
t(x,x,x) := 6 \det(x) \tbtmat{-\alpha}{j}{-j'}{\beta}.
\]
Then $(V,b,t)$ is a Freudenthal triple system.  Since
\[
q(x,x,x,x) := 12 \det(x)^2,
\]
it is certainly degenerate, and we denote it by $\M_s$.
By \cite[\S 4]{Brown:E7} or \cite[\S 4]{Meyb:FT}, all degenerate
triple systems are forms of one of these, meaning that any degenerate
triple system becomes isomorphic to $\M_s$ over a separable closure of
the ground field.
\end{eg}

\begin{eg} \label{nondegeg} 
For $J$ an Albert $F$-algebra, there is a nondegenerate triple system
denoted by $\M(J)$ whose underlying $F$-vector space is $V =
\stbtmat{F}{J}{J}{F}$.  Explicit formulas for the $b$, $t$, and $q$
for this triple system can be found in \cite[\S 3]{Brown:E7}, \cite[\S 
6]{Meyb:FT}, \cite[\S 1]{Ferr:strict}, and \cite[3.2]{G:struct}.  For
$J^d$ the split Albert $F$-algebra, we set $\Md := \M(J^d)$.  It is
called the {\em split} triple system because $\Inv(\Md)$ is the split
simply connected  
algebraic group of type $E_7$ over $F$ \cite[3.5]{G:struct}.  By
\cite[\S 4]{Brown:E7} or \cite[\S 4]{Meyb:FT} every nondegenerate
triple system is a form of $\Md$.
\end{eg}

A {\em similarity} of triple systems is 
a map $f \!: (V, b, t) \iso (V', b', t')$ defined by an $F$-vector
space isomorphism $f \!: V \iso V'$ such that $b'(f(x), f(y)) =
\lambda b(x,y)$ and $t'(f(x), f(y), f(z)) = \lambda f(t(x,y,z))$ for
all $x$, $y$, $z \in V$ and some $\lambda \in \Fx$ called the {\em
multiplier} of $f$.  Similarities with multiplier one are called {\em
isometries}. They are the isomorphisms in the obvious category of
Freudenthal triple systems.
For a
triple system $\M$, we
write $\Inv(\M)$ for the algebraic group whose $F$-points are the
isometries of $\M$.

\begin{rmk} \label{degender}
Although it is not clear precisely what the structure of the
automorphism group of a degenerate triple system is, a few simple 
observations can be made which make it appear to be not very
interesting from the standpoint of simple algebraic groups.

Since by definition any element of $\Inv(\Ms)$ must preserve the
quartic form $q$, it must also be a similarity of the quadratic form
$\det$ with multiplier $\pm 1$.  This defines a map $\Inv(\M_s)
\longto \mmu2$ which is surjective since $\varpi \in \Inv(\M_s)(F)$
maps to $-1$, where
\[
\varpi \basmat := \tbtmat{-\beta}{j'}{j}{\alpha}.
\]
So $\Inv(\Ms)$ is not connected.

Also, we can make some bounds on the dimension.  Specifically, we
define a map $f \!: \GmF \times W \times GL(W) \longto \Inv(\Ms)$ by
\[
f(c, u, \phi) \basmat := \tbtmat{c \alpha}{\phi(j)}{\alpha u +
  \phi^\dagger(j')}{\frac1c (\beta + s(\phi(j), u))},
\]
where $\phi^\dagger = \s(\phi)^{-1}$ for $\s$ the involution on
$\EndF(W)$ which is adjoint with respect to $s$.  (So $s(\phi(w),
\phi^\dagger(w')) = s(w, w')$ for all $w, w' \in W$.)  Then $f$ is an
injection of varieties, but 
\[
f(c, u, \phi) f(d, v, \psi) = f(c d, du+\phi^\dagger(v), \phi \psi),
\]
so it is not a group homomorphism.
It does, however, restrict to be a morphism of algebraic groups on
$\GmF \times \{0\} \times GL(W)$, so $\Inv(\M_s)$ contains a split
torus of rank 28.
This map $f$ is also not surjective since for any $u \ne 0$, the map
\[
\varpi f(c, u, \phi) \varpi^{-1} \basmat = 
\tbtmat{\frac1c (\alpha - s(\phi(j'), u))}{\beta u + \phi^\dagger(j)}%
{\phi(j')}{\beta c}
\]
is not in the image of $f$.
For an upper bound, we observe that the
identity component of $\Inv(\Ms)$ is contained in the isometry group
of $\det$ (whose identity component is of type $D_{28}$, hence is 
1540-dimensional).
Thus
\[
757 < \dim \Inv(\M_s)^{+} \le 1540.
\]
\end{rmk}

\section{An identity} \label{pisec}

For a Freudenthal triple system $(V, b, t)$ over $F$, we define an $F$-vector
space map $p \!: V \otimes_F V \longto \EndF(V)$ given by
\begin{neqn} \label{pdef}
p(u \otimes v) w := t(u,v,w) - b(w,u)v - b(w,v)u.
\end{neqn}
In the case where the triple system is nondegenerate, Freudenthal
\cite[4.2]{Frd:E7.1} also defined a map $V \otimes_F V 
\longto \EndF(V)$ which he denoted by $\times$.  The obvious
computation shows that his map is related to our map $p$ by
\begin{neqn} \label{Frdrelate}
8\,v \times v' = p(v \otimes v').
\end{neqn}

\begin{thm} \label{pithm}
Let $\M := (V, b, t)$ be a Freudenthal triple system with map $p$ as
given above.   Then $\M$ is nondegenerate if and only if it satisfies
the identity
\begin{neqn} \label{pi}
\tr (p(x \otimes x)\, p(y \otimes y)) = 24 \,\big(q(x,x,y,y) - 2\, b(y,x)^2\big)
\end{neqn}
for all $x, y \in V$, where $\tr$ is the usual trace form on $\EndF(V)$.
\end{thm}

\begin{pf}
If $\M$ is nondegenerate, then the conclusion is 
\cite[31.3.1]{Frd:E7.10} or it can be easily derived from
\cite[7.1]{Meyb:FT}.  So we may assume that $\M$ is degenerate and
show that it doesn't satisfy (\ref{pi}).  Extending scalars, we may
further assume that our ground field is separably closed and so that $\M =
\M_s$, the degenerate triple system from Example \ref{degeneg}.

To simplify some of our
formulas, we define the {\em weighted determinant}, $\wdet \!: \M_s
\longto F$, to be given by
\[
\wdet(x) := 3\alpha \beta - s(j, j')
\]
for $x$ and $y$ as in (\ref{xydef}).

We first compute the value of the left side of (\ref{pi}).  For $x$
and $y$ as in (\ref{xydef}), we can
directly calculate the action of $p(x \otimes x)\, p(y \otimes y)$ on 
each of the four entries of our matrix as in (\ref{matdef}).  Since we 
are interested in the trace of this operator, we only record the
projection onto the entry that we are looking at.
\begin{neqn} \label{ultr}
\tbtmat{1}{0}{0}{0} \mapsto 4 \left[ \wdet(x)
  \wdet(y) - 4 \alpha \delta s(j, k') \right]
\end{neqn}
\begin{neqn} \label{lrtr}
\tbtmat{0}{0}{0}{1} \mapsto 4 \left[ \wdet(x) \wdet(y) + 4 \beta
  \gamma s(j', k) \right]
\end{neqn}
\begin{neqn}\label{urtr}
\tbtmat{0}{m}{0}{0} \mapsto 4 \left[ 
\begin{array}{c}
\det(x) \det(y) m + 4 (\alpha
  \delta - s(k, j') )\, s(k', m) j \\
+ 2 \det(x) s(k', m) k + 2 \det(y)
  s(j', m) j 
\end{array}
\right] 
\end{neqn}
\begin{neqn} \label{lltr}
\tbtmat{0}{0}{m'}{0} \mapsto 4 \left[ 
\begin{array}{c}
\det(x) \det(y) m' - 4 (\beta
  \gamma - s(j, k'))\, s(k, m') j' \\
- 2 \det(x) s(k, m') k' - 2 \det(y)
    s(j, m') j'
\end{array}
\right]
\end{neqn}

Since $s$ is nondegenerate, it induces an identification of $V$ with
its dual $V^\ast$ by sending $x \in V$ to the map $v \mapsto s(x,
v)$.  We may also identify $V \otimes_F V^\ast$ with $\EndF(V)$, and
combining these two identifications provides an isomorphism $\varphi_s 
\!: V \otimes_F V \iso \EndF(V)$ given by
\begin{neqn} \label{phib}
\varphi_s(x \otimes y) w = x s(y, w),
\end{neqn}
cf.~\cite[5.1]{KMRT}.  One has $\tr(\varphi_s(x \otimes y)) = s(y, x)$.

With that notation, the terms in the brackets of (\ref{urtr}) and
(\ref{lltr}) with coefficient 4 give the maps
\[
(\alpha \delta - s(k, j')) \varphi_s(j \otimes k') \hbox{\ and\ }
-(\beta \gamma - s(j, k')) \varphi_s(j' \otimes k)
\]
on $W$,
which have traces
\[
(\alpha \delta -s(k, j')) s(k', j) \hbox{\ and\ } 
(\beta \gamma - s(j, k')) s(j', k)
\]
respectively.
Similarly, the terms with coefficient 2 give the maps
\[
\det(x)\varphi_s(k \otimes k') + \det(y) \varphi_s (j \otimes j')
\hbox{\ and\ } 
-\det(x) \varphi_s(k' \otimes k) - \det(y) \varphi_s (j' \otimes j)
\]
whose sum has trace
\[
2 \left[ \det(x) s(k',k) + \det(y) s(j', j) \right].
\]
Adding all of this up, we see that
\begin{neqn} \label{trform}
\begin{array}{rcl}
\frac18 \tr(p(x \otimes x) p(y \otimes y)) &=& \wdet(x) \wdet(y) + 2
\beta \gamma s(j', k) - 2 \alpha \delta s(j, k') \\
&& +\, 27 \det(x) \det(y)
\\
&& +\, 2 \left[ (\beta \gamma - s(j, k')) s(j', k) + (\alpha \delta -
  s(k, j')) s(k', j) \right] \\
&& +\, 2 \left[ \det(x) s(k', k) + \det(y) s(j', j) \right] \\
&=& 3 q(x,x,y,y) - b(y,x)^2 + 20 \det(x) \det(y) - 5 \det(x,y)^2,
\end{array}
\end{neqn}
where we have linearized the determinant to define
\[
\det(x, y) := \det(x + y) - \det(x) - \det(y).
\]

Taking the difference of the last expression in (\ref{trform}) and
one-eighth of the right side of (\ref{pi}), we get the quartic
polynomial 
\begin{neqn} \label{remainder}
5\,b(x,y)^2 + 20 \det(x) \det(y) - 5\det(x,y)^2.
\end{neqn}
We plug
\[
x := \tbtmat{0}{j}{j'}{0} \hbox{\ and\ } y := \tbtmat{0}{k}{k'}{0}
\]
into (\ref{remainder}),
where we have chosen $j$, $j'$, $k$, and $k'$ such that
\[
s(j, k') = s(j', k) = 0 \hbox{\ and\ } s(j, j') = s(k, k') = 1.
\]
Then $b(x,y) = 0$ and $\det(x) = \det(y) = -1$.  Since
$\det(x + y) = -2$, we have $\det(x,y) = 0$.  Thus the value of
(\ref{remainder}) is 20 and (\ref{pi}) does not hold for degenerate
triple systems.
\end{pf}

It really was important that we allow $x \ne y$ in
(\ref{pi}), since for $x = y$ (\ref{remainder}) is identically zero.  
So all triple systems satisfy
the identity
\begin{neqn} \label{badtrid}
\tr (p(x \otimes x)^2) = 24 \, q(x,x,x,x).
\end{neqn}

\section{Gifts} \label{giftsec}

In this section we will define an object we call a gift, such that
every adjoint group  
of type $E_7$ is the automorphism group of some gift.  We must first
have some preliminary definitions.

Suppose that $\As$ is a central simple algebra with a symplectic
involution $\s$.  The sandwich map
\[
\Sand \!: A \otimes_F A \longto \EndF(A)
\]
defined by 
\[
\Sand (a \otimes b) (x) = axb \hbox{\, for $a, b, x \in A$}
\]
is an isomorphism of $F$-vector spaces by \cite[3.4]{KMRT}.  Following
\cite[\S 8.B]{KMRT}, we define a map $\st$ on $A \otimes_F A$
which is defined implicitly by the equation
\[
\Sand (\st(u))(x) = \Sand(u)(\s(x)) \hbox{\ for $u \in A \otimes A$,
$x \in A$}.
\]

Suppose now that $A$ is split.  Then $A \cong \EndF(V)$ for some
$F$-vector space $V$, and $\s$ is the adjoint involution for some
nondegenerate skew-symmetric bilinear form $b$ on $V$ (i.e., $b(fx, y)
= b(x, \s(f)y)$ for all  
$f \in \EndF(V)$).
As in (\ref{phib}), we have an identification $\varphi_b \!: V
\otimes_F V \iso \EndF(V)$, and by a straightforward computation (or
see \cite[8.6]{KMRT}), $\st$ 
is given by 
\begin{neqn} \label{st}
\st(\varphi_b(x_1 \otimes x_2) \otimes \varphi_b(x_3 \otimes x_4)) =
- \varphi_b(x_1 \otimes x_3) \otimes \varphi_b(x_2 \otimes x_4)
\end{neqn}
for $x_1$, $x_2$, $x_3$, $x_4 \in V$.  

Finally, for $f \!: A \longto A$ an $F$-linear map, we define
$\widehat{f} \!: A \otimes_F A \longto A$ by 
\[
\widehat{f}(a \otimes b) = f(a)b.
\]

\begin{defn} \label{giftdef}
A {\em gift} $\G$ over a field $F$ is a triple $\bgift$ such that $A$ is a
central simple $F$-algebra of degree 56, $\s$ is a symplectic
involution on $A$, and $\pi \!: A \longto A$ is an $F$-vector space
map such that
\begin{description}
\item[G1] $\s \pi(a) = \pi \s(a) = -\pi(a)$,
\item[G2] $a \pi(a) \ne 2a^2$ for some $a \in \Skew\As$, 
\item[G3] $\pi(\pi(a) a) = 0$ for all $a \in \Skew\As$, 
\item[G4] $\pihat - \shat - \idhat = - (\pihat - \shat - \idhat) \s_2$, and
\item[G5] $\TrdA(\pi(a)\pi(a')) = -24 \TrdA(\pi(a) a')$ for all $a, a'
  \in \As$. 
\end{description}
\end{defn}

By $\Skew\As$ we mean the vector space of $\s$-skew-symmetric elements
of $A$, i.e., those $a \in A$ such that $\s(a) = -a$.

A definition as strange as \ref{giftdef} demands an example.

\begin{eg} \label{endeg}
Suppose that $\M = (V,b,t)$ is a nondegenerate Freudenthal triple
system over $F$.  
Set $\EndM := (\EndF(V), \s, \pi)$ where $\s$ is the involution
on $\EndF(V)$ adjoint to $b$.  Using the identification $\varphi_b \!: V \otimes
V \longto \EndF(V)$, we define $\pi \!: A \longto A$ by $\pi := p
\varphi^{-1}_b$, where $p$ is as in (\ref{pdef}). 

We show that $\EndM$ is a gift.
A quick computation shows that 
\[
-b(\pi(\varphi_b(x \otimes y)) z, w) = b(z, \pi(\varphi_b(x \otimes
y)) w) = -b(\pi(\varphi_b(y \otimes x)) z, w), 
\]
which demonstrates G1, since $\s\varphi_b(x \otimes y) = -\varphi_b(y
\otimes x)$. 

Suppose that G2 fails.  Then for $v \in V$, we set $a := \varphi_b(v
\otimes v)$ and observe that $a^2 = 0$, so
\[
0 = \varphi_b(v \otimes v) \pi (\varphi_b(v \otimes v)) v 
= q(v,v,v,v) v.
\]
Since this holds for all $v \in V$, $q$ is identically
zero, contradicting FTS2.  Thus G2 holds.

Since elements $a$ like in the preceding paragraph span $\Skew\As$, in
order to prove G3 we may show that
\begin{neqn} \label{G3eqn}
\pi \left(
\pi( a ) a' + \pi( a' ) a
\right) y = 0,
\end{neqn}
where 
\begin{neqn} \label{aaprime}
a = \varphi_b(x \otimes x) \hbox{\ and\ } a' = \varphi_b(z \otimes z).
\end{neqn}
A direct expansion of the left-hand side of (\ref{G3eqn})
shows that it is equivalent to FTS3$^\prime$.

Using just the bilinearity and skew-symmetry of $b$ and the
trilinearity of $t$, G4 is equivalent to
\[
\hbox{$t(x,y,x') = t(x,x',y)$ for all $x, x', y \in V$.}
\]
Thus G4 holds by FTS1.

Finally, consider G5.  If $a$ is symmetric, then by G1 $\pi(a) = 0$
and the identity holds.  If $a'$ is symmetric then the left-hand side
of G5 is again zero by G1.  Since $\s$
and $\TrdA$ commute, we have
\[
\TrdA(\pi(a) a') = \s (\TrdA(\pi(a)a')) = -\TrdA(a' \pi(a)) =
-\TrdA(\pi(a) a'),
\]
so the right-hand side of G5 is also zero.  Consequently, by the
bilinearity of G5, we may assume that $a$ and $a'$ are   
skew-symmetric, and we may further assume that
$a$ and $a'$ are as given in (\ref{aaprime}).  Then G5 reduces to
(\ref{pi}). 
\end{eg}

It turns out that this construction produces all
gifts
with split central simple algebra component.

\begin{lem} \label{endlem}
Suppose that $\G = \bgift$ is a gift over $F$.  Then $\G \cong \EndM$
for some nondegenerate Freudenthal triple system over $F$ if and only
if $A$ is split. 
\end{lem}

\begin{pf}
One direction is done by Example \ref{endeg}, so suppose that $\bgift$
is a gift with $A$ split.  Then we 
may write $A \cong \EndF(V)$ for some 56-dimensional $F$-vector space
$V$ such that $V$ is endowed with a nondegenerate skew-symmetric
bilinear form $b$ and $\s$ is the involution on $A$ which is adjoint
to $b$.  We define $t \!: V \times V \times V \longto V$ by
\[
t(x,y,w) := \pi(\varphi_b(x \otimes y)) w + b(w,x) y + b(w,y)x.
\]
Observe that $t$ is trilinear.  We define a 4-linear form $q$ on $V$
as in FTS2.

The proof that FTS3$^\prime$ implies G3 in Example \ref{endeg} reverses to
show that G3 implies FTS3$^\prime$.  Similarly, G4 implies that $t(x,y,z) =
t(x,z,y)$ for all $x,y,x' \in V$, so $q(w,x,y,z) = q(w,x,z,y)$.  G1
implies that $q(w,x,y,z) = q(z,x,y,w) = q(w,y,x,z)$.  Since the
transpositions
$(3\,4)$, $(1\,4)$, and $(2\,3)$ generate
$\mathcal{S}_4$ (= the symmetric group on four letters) $q$ is 
symmetric.

Next, suppose that FTS2 fails, so that
$q$ is identically zero.  Then since $b$ is nondegenerate, $t$ is also
zero.  Then for $v, v', z \in V$ and $a := \varphi_b(v \otimes v)$ and
$a' := \varphi_b(v' \otimes v')$, 
\[
(a \pi(a') + a' \pi(a))z = 2 (b(v,v') b(v',z) v + b(v',v) b(v,z)) =
2(aa' + a' a) z.  
\]
Since elements of the same form as $a$ and $a'$ span $\Skew\As$, this
implies that G2 fails, which is a contradiction.  Thus FTS2 holds and
$(V,b,t)$ is a Freudenthal triple system.

Finally, writing out G5 in terms of $V$ gives (\ref{pi}), which shows
that $(V,b,t)$ is nondegenerate.
\end{pf}

\begin{rmk}
The astute reader will have noticed that our definition of $\EndM$
almost works if $\M$ is degenerate, in that the only problem is that the
resulting $(A, \s, \pi)$ doesn't satisfy G5.  That example and the
proof of \ref{endlem} make it clear that if we remove the axiom G5
from the definition of a gift, then we would get an analog to Lemma
\ref{endlem} where the Freudenthal triple system is possibly degenerate.
\end{rmk}

\begin{rmk} \label{simrmk}
Observe that in the isomorphism $\G \cong \EndM$ from the preceding 
lemma, $\M$ is only determined up to similarity.  To wit,
for a nondegenerate triple system $\M = (V, b, t)$ and $\lambda \in
\Fx$, we define  
a similar   
structure $\M_\lambda = (V, \lambda b, \lambda t)$.  Then
$\M_\lambda$ is also a nondegenerate triple system and $\EndM =
\End{}(\M_\lambda)$.  The only  
potential difficulty with this last equality would be if 
the $\pi$ produced by $\M_\lambda$, which we
shall denote by $\pi_\lambda$, is different from the $\pi$ produced by 
$\M$.  However, we see that
\[
\begin{array}{rcl}
\pi_\lambda (\varphi_{\lambda b} (x \otimes y)) w &=& \lambda t(x,y,w) -
\lambda b(w,x) y -
\lambda b(w,y)x \\
&=& \lambda \pi (\varphi_{b} (x \otimes y)) w \\
&=& \pi (\varphi_{\lambda b} (x \otimes y)) w.
\end{array}
\]
\end{rmk}

\subsection*{Isometries and derivations}

\begin{defn}
An {\em isometry} in a gift $\G := \bgift$ is an element $f \in A$
such that $\s(f) f = 1$ and $\pi(f a f^{-1}) = f \pi(a) f^{-1}$ for all $a \in
A$ (this ensures that $\Int(f)$ is an automorphism of $\G$).  The
motivation for the name is that then the isometries in $\EndM$ are
just the isometries of $\M$.
We
set $\Iso(\G)$ to  
be the algebraic group whose $F$-points are the group of isometries in
$\G$.  

A {\em derivation} in $\G$ is an element $f \in \Skew\As$ which satisfy
\begin{description}
\item[GD] $\pi(fa) - \pi(a f) = f \pi(a) - \pi(a) f$ for all $a \in A$.
\end{description}
We define $\Der(\G)$ to be the vector space of derivations in $\G$.
The name derivation has the same sort of motivation: We will see in
the proof of Proposition \ref{impi} that the derivations
in $\End{}(V,b,t)$ are precisely the maps in $\EndF(V)$ which are  
traditionally called derivations of the triple system $(V,b,t)$.
\end{defn}

\begin{borel*}
The description of the isometries in $\End{}(V,b,t)$ combined
with Lemma \ref{endlem} shows that $\Iso(\G)$ is simple simply
connected of type $E_7$.
Since any automorphism of $\G$ is also an isomorphism of $A$, it must
be of the form $\Int(f)$ for some $f \in \Ax$ such that $\s(f)f \in
\Fx$.
Thus the map $\Iso(\G)
\longto \aut(\G)$ is a surjection over a separable closure of the
ground field.  Since the kernel of
this map is the center of $\Iso(\G)$, $\aut(\G)$ is simple adjoint of
type $E_7$.

We can actually say more.  The group $\Iso(\G)$ is a subgroup of the
symplectic 
group $\Sp\As$, whose $F$-points are the elements $f \in A$ such
that $\s(f) f = 1$.  It is this embedding combined with the vector
representation of $\Sp\As$ (= the natural embedding of $\Sp\As$ in
$A^\ast$) which gives rise to the Tits algebra of
$\Iso(\G)$, and so the Tits algebra of $\Iso(\G)$ is the same as the
Tits algebra of $\Sp\As$, which is just $A$.

It is easy to see that $\Der(\G)$ is actually a Lie subalgebra of
$\Skew\As$, where the bracket is the usual commutator.  In fact, by
formal differentiation as in \cite[3.21]{Borel} or 
\cite[\S 4]{Jac:J1}, $\Der(\G)$ is the Lie algebra of $\Iso(\G)$. 
\end{borel*}

\begin{prop} \label{impi}
For $\G := \bgift$ a gift, $\im \pi = \Der(\G)$.
\end{prop}

\begin{pf}
Since $\im \pi$ and $\Der(\G)$ are both vector subspaces of $A$, it is 
equivalent to prove this over a separable closure.  Thus we may assume 
that $A$ is split, so that $\G = \EndM$ for some nondegenerate
triple system $\M := (V,b,t)$ over $F$ by Lemma \ref{endlem}.

Consider the vector subspace $D$ of $\Skew\As$ consisting of elements
$d$ such that 
\[
dt(u,v,w) = t(du, v, w) + t(u, dv, w) + t(u, v, dw)
\]
for all $u, v, w \in V$.  (These are
known as the derivations of $\M$.)  The obvious computation shows that
$\im \pi \subseteq D$, which one can find in \cite[p.~166,
Lem.~1.3]{Meyb:FT}.  Conversely, the reverse containment holds by
\cite[p.~185, S.~8.3]{Meyb:FT}.  (He has an ``extra'' hypothesis that
the characteristic of $F$ is not 5 because he is also considering
triple systems of dimensions 14 and 32, but that is irrelevant for our 
purposes.)

For $d \in \Skew\As$, consider the element
\begin{neqn} \label{adid}
\pi(da) - \pi(ad) - d\pi(a) + \pi(a) d 
\end{neqn}
in $A$, which is zero if and only if $d$ satisfies GD.  Since $\pi$ is
linear, \ref{adid} is zero for all $a$ if and only if it is zero when
$a = \varphi_b(u \otimes v)$.  Applying the endomorphism
\ref{adid} to $w \in V$ and expanding out, we obtain
\[
t(du, v, w) + t(dv, u, w) + t(u,v,dw) - dt(u,v,w).
\]
So $d$ satisfies GD if and only if it lies in $D = \im \pi$.
\end{pf}

\subsection*{A category equivalence}

We will now show that there is an equivalence of
categories between the category of adjoint groups of type
$E_7$ over $F$ and the category of gifts over $F$, where both
categories have isomorphisms for morphisms (i.e., they are groupoids).
We use the notation and vocabulary of \cite[\S 26]{KMRT} with impunity.

\begin{thm} 
The automorphism group of a gift defined over
a field $F$ is an adjoint group of type $E_7$ over $F$.  This provides
an equivalence between the groupoid of gifts over $F$ and the groupoid 
of adjoint groups of type $E_7$ over $F$.
\end{thm}

\begin{pf}
Let $\C_\G(F)$ denote the groupoid of gifts over $F$
and let $\C_{E_7}(F)$ denote the groupoid of adjoint groups of type
$E_7$ over $F$.  Let
\[
S(F) \!: \C_\G(F) \longto \C_{E_7}(F)
\]
be the functor induced by the map on objects given by $\G \mapsto
\aut(\G)$.  Then we have a commutative diagram
\[
\begin{CD}
\C_\G(F) @>{S(F)}>> \C_{E_7}(F)\\
@V{i}VV @VV{j}V \\
\C_\G(F_s) @>{S(F_s)}>> \C_{E_7}(F_s) 
\end{CD}
\]
where $F_s$ is a fixed separable closure of $F$ and $i$, $j$ are the
obvious scalar extension maps.  These maps are both 
``$\Gamma$-embeddings'' for $\Gamma$ the Galois group of $F_s$ over
$F$, in that there is a $\Gamma$-action on the morphisms in the
category over $F_s$ with fixed points the morphisms coming from $F$.
Since the diagram is commutative and is compatible with the $\Gamma$-action on
the morphisms in the categories over $F_s$, $S(F_s)$ is said to be a
``$\Gamma$-extension'' of $S(F)$.

Since any nondegenerate triple system is a form of $\Md$
(\ref{nondegeg}), every nondegenerate gift is a form of $\End{}(\M^d)$
by Lemma \ref{endlem}.  Thus $\C_\G(F_s)$ is
connected.   Since
$\C_{E_7}(F_s)$ is also connected and any object in either category has
automorphism group the split adjoint group of type $E_7$, $S(F_s)$ is
an equivalence of groupoids.  (For $G$ adjoint of type $E_7$,
$\Aut{}(G) \cong G$ by \cite[1.5.6]{Ti:Cl}.)

By \cite[26.2]{KMRT} we only need to show that $i$ {\em satisfies the descent
condition}, i.e., that $1$-cocycles in the automorphism group of some
fixed element of $\C_\G(F_s)$ define objects 
in $\C_\G(F)$.
Let $(A, \s, \pi)$ be a gift over $F$ and set
\[
W := \HomF(A \otimes_F A, A) \oplus \HomF(A,A) \oplus \HomF(A,
\Skew\As).
\]
Here $(m, \s, \pi) \in W$ gives the structure of $(A, \s, \pi)$.  
The rest of the argument is as in 26.9, 26.12, 26.14, 26.15, 26.18, or 
26.19 of \cite{KMRT}.  We let $\rho$ denote the natural map
$GL(A)(F_s) \longto GL(W)(F_s)$ and observe that elements of the orbit 
of $w$ under $\im \rho(F_s)$ define all possible gifts over $F_s$ and
the objects of $\C_\G(F)$ can be identified with the set of all $w'
\in W$ such that $w'$ is in the orbit of $w$ over $F_s$.  Then $i$
satisfies the descent condition by \cite[26.4]{KMRT}.
\end{pf}

\section{Applications to flag varieties} \label{flagsec}
\renewcommand{\C}{\mathfrak{C}}

\begin{defn} \label{idealdef}
An {\em inner ideal} of a gift $\G = \bgift$ is a right ideal $I$ of
$A$ such that $\pi(I\s(I)) \subseteq I$.  A {\em singular ideal} of
$\G$ is a right ideal of $A$ such that $\pi(I\s(I)) = 0$.
\end{defn}

If $A$ is split so that $\G \cong \EndM$ for some triple system $\M =
(V,b,t)$, there is a bijection between subspaces $U$ of $V$ and right
ideals $\Hom_F(V,U)$ of $\EndF(V)$.  In this bijection, inner ideals
of $\M$ (i.e., those subspaces $U$ such that $t(U,U,V) \subseteq U$)
correspond to inner ideals of $\G$.  The same statement holds for singular
ideals where a singular ideal of $\M$ is defined to be a subspace $U$
such that $t(u, u', v) = b(v, u)u' + b(v, u')u$ for all $u, u' \in U$
and $v \in V$.

Since all inner ideals in a nondegenerate Freudenthal triple system
are totally isotropic with respect to the skew-symmetric bilinear form 
\cite[2.4]{Ferr:strict}, any singular or inner ideal in a gift $\G$ is 
isotropic, i.e., $\s(I)I = 0$.

Now all of the flag varieties (a.k.a.~homogeneous projective
varieties) for an arbitrary group $\aut(\G)$ of type $E_7$ can easily
be described in terms of the singular and rank 12 inner ideals of the
gift $\G$, by translating the corresponding description for the flag
varieties for such groups with trivial Tits algebras from
\cite[\S 7]{G:struct}.  Specifically, one simply takes the statement
of \cite[7.5]{G:struct} and replaces every instance of
``$n$-dimensional'' with ``rank $n$'' as well as replacing
$\Inv(\B)$ with $\aut(\G)$.

\section{A construction} \label{constsec}

The rest of this section is taken up with a construction which
produces groups of type $E_7$ with Tits algebras of index 2, and in
particular all gifts (hence all groups of type $E_7$) over a
real-closed field.  (We say that a field $F$ is {\em real} if
$-1$ is not a sum of squares in $F$ and that it is {\em real-closed}
if it is real and no algebraic extension is real, please see
\cite[Ch.~9, \S\S 1, 2]{Lam} for 
more information.)  Since one of these groups is anisotropic, this is
the first construction of an anisotropic group of type $E_7$ directly
in terms of this 56-dimensional form. (Groups of this type
have been implicitly constructed as the  
automorphism groups of Lie algebras given by the Tits construction.)  

\subsection*{A diversion to symplectic involutions}

Let $Q := (\alpha, \beta)_F$ be the quaternion algebra over $F$ generated by
elements $i$ and $j$ such that  $i^2 = \alpha$, $j^2 = \beta$, and $ij = -ji$.
Our
construction will involve taking a 56-dimensional skew-symmetric
bilinear form from a triple system over $F$ and twisting it to get
a symplectic involution $\s$ on $M_{28}(Q)$.  However, if our triple system
is of the form $\M(J)$ for some $J$, then the skew-symmetric form is of a very
special kind, namely it is obtainable from a quadratic form
$q$ (specifically, $\qform{1} \perp T$ where $T$ is the trace on the Jordan 
algebra), and we use this fact to get an explicit description of $\s$.

More generally, suppose that $(V,q)$ is a nondegenerate quadratic
space over $F$ with associated symmetric bilinear form $b_q$ such that 
$b_q(x,x) = q(x)$.  Let $K = F(\sqrt{\alpha})$ with nontrivial
$F$-automorphism $\iota$ and fix some square root of $\alpha$ in $K$.
Then for $W_K = (V \oplus V) \otimes K$, we set $s$ to be
the skew-symmetric bilinear form given by 
\[
s(w_1, w_2) := \sqrt{\alpha} \left[ b_q(v_1, v'_2) - b_q(v_2, v'_1) \right],
\]
where $w_i = (v_i, v'_i) \in W_K$ for $i = 1, 2$.
Since $q$ is nondegenerate, $s$ is as well.

We define a twisted ction of $\iota$ on $W_K$ by setting $^\iota w = (\psi_c
\otimes \iota) w$ where the $\iota$ on the right denotes the usual
action and $\psi_c$ is given by $\psi_c(v, v') = (cv', (\beta/c) v)$
for some $c \in \Fx$.  Then the fixed subspace is a vector space over
$F$ which we denote by $W$.  Note that $W \otimes_F K \cong W_K$.
Since $s(^\iota w_1, ^\iota w_2) = \beta \iota s(w_1, w_2)$, the
action of $\iota$ on $W_K$ induces an $\iota$-semilinear automorphism
$\Int(\psi_c) \otimes \iota$ of $(\EndK(W_K), \s_s)$.

\begin{lem} \label{symplem}
Suppose that for $1 \le i \le n$, $(V_i,q_i)$ is a nondegenerate
quadratic space over $F$ isomorphic to $\qform{a_i}$, $K =
F(\sqrt{a})$,
and  $(W_{K,i}, s_i)$ is the associated symplectic
space defined over $K$ as described above, and for $c_i \in \Fx$, $(A, \s)$ is
fixed subalgebra defined by descent from $(\EndK(W_K), \s_s)$ by
$\Int(\psi_{c_1} \osum \cdots \osum \psi_{c_n}) \otimes \iota$.  Then $A$ is
Brauer-equivalent to $Q = (\alpha, \beta)_F$ and 
the involution $\s$ is adjoint to the hermitian form 
\[
\qform{c_1 a_1, \ldots, c_n a_n}
\]
over $Q$.
\end{lem}

There is a unique symplectic involution on $Q$ which we will denote by
$\gamma$.  It is defined by setting $\gamma(i) = -i$ and $\gamma(j) = -j$.
In the statement of the lemma, {\em hermitian} 
means with respect to the involution $\gamma$.

\medskip
\begin{pf}
Consider the two isomorphisms
\[
(Q \otimes_F K) \otimes_F M_n(F) \xrightarrow{f} M_2(K) \otimes_F
M_n(F)
\xrightarrow{g} M_{2n}(F) \otimes_F K.
\]
We fix a square root of $\alpha$ in $K$ and 
let $E_{rs}$ be the
matrix whose only nonzero entry is a 1 in the $(r,s)$-position. 
We define $f$ by  
\[
f(1 \otimes E_{rs}) := \tbtmat{1}{}{}{\frac{c_s}{c_r}} \otimes E_{rs},
\]
\[
f(i \otimes 1) :=
\tbtmat{\sqrt{\alpha}}{0}{0}{-\frac{c_r}{c_s}\sqrt{\alpha}} \otimes 1,
\]
and
\[
f(j \otimes E_{rs}) := \tbtmat{0}{c_s}{\beta/c_r}{0} \otimes
E_{rs}.
\]
We define $g$ by dividing $M_{2n}(K)$ into $n\times n$ blocks and
setting $g(x \otimes E_{rs})$ to be the matrix with $x$ in the $(r,s)$
block.

Set $m$ in $M_{2n}(F) \otimes K$ to be the matrix with diagonal blocks
$C_1, \ldots, C_n$ for $C_i := \stbtmat{0}{c_i}{\beta/c_i}{0}$.  Then $A$
is the $F$-subalgebra of $M_{2n}(F) \otimes_F K$ fixed by $\Int(m)
\otimes \iota$, for $\iota$ the nontrivial $F$-automorphism of $K$.
Since $g^{-1}(\Int(m) \otimes \iota)g$ fixes $f(Q \otimes F
\otimes M_n(F))$, $gf$ provides an $F$-isomorphism $Q \otimes_F F \otimes_F
M_n(F) \cong A$.

The involution $\tau$ on $M_{2n}(F) \otimes K$ which is adjoint relative to
$s_1 \oplus \cdots \oplus s_n$ is 
\[
\tau := \Int\thbthmat{A_1^{-1}}{}{}{}{\ddots}{}{}{}{A_n^{-1}} \circ
\Int\thbthmat{J}{}{}{}{\ddots}{}{}{}{J} \circ t,
\]
where $A_i = \sqrt{\alpha} a_i \Id_2$, $J = \stbtmat{0}{1}{-1}{0}$, and $t$ is
the transpose \cite[p.~24]{KMRT}.  For $\gamma$ the unique
symplectic involution on $Q$, we have an involution on $M_{2n}(F)
\otimes K$ given by
\[
(gf)(\gamma \otimes \Id_K \otimes t)(gf)^{-1} =
\left(\Int\thbthmat{c_1^{-1} J}{}{}{}{\ddots}{}{}{}{c_n^{-1} J} \circ t\right) \otimes \Id_K. 
\]
So the involution $\tau$ induced by $\s$ on $Q \otimes K \otimes
M_n(F)$ satisfies
\[
(gf)\s(\gamma \otimes \Id_K \otimes t)(gf)^{-1} = \tau (gf) (\gamma
\otimes \Id_K \otimes t) (gf)^{-1} = \Int \sthbthmat{c_1
A_1^{-1}}{}{}{}{\ddots}{}{}{}{c_n A_n^{-1}}.
\]
Since $(gf)^{-1} (c_i A_i^{-1}) = c_i (\sqrt{\alpha} a_i)^{-1}$, the
lemma is proven.
\end{pf}

\subsection*{The construction}
In characteristic $\ne 5$ (and as always $\ne 2, 3$), we can define a
{\em Brown algebra} to be a 
56-dimensional central simple structurable 
algebra with involution $\bbrown$  such that the
space of skew-symmetric elements is one-dimensional.  In
characteristic 5, a different definition is needed at the moment due
to insufficiently strong classification results in that
characteristic.  See \cite[\S 2]{G:struct} for a full definition and
\ref{browneg} below for examples. 

The relevant point is that given a Brown $F$-algebra $\B$, one can produce
a nondegenerate Freudenthal triple system $\M := (V, b, t)$ in a
relatively natural way, see \cite[\S 2]{AF:CD} or \cite[\S
4]{G:struct}.  This triple system is determined only up to similarity
(i.e., for every $\lambda \in \Fx$, $(V, \lambda b, \lambda t)$ is
also a triple system associated to $\B$, and these are all of them).
Also, this construction produces {\em all} nondegenerate Freudenthal
triple systems over $F$ by \cite[4.14]{G:struct}.  

We define a gift
$\EndB$ by setting $\EndB := \EndM$.  Although $\M$ is only determined 
up to similarity by $\B$, $\EndB$ is still well-defined, as observed in 
Remark \ref{simrmk}.

Our construction will need to use a specific kind of Brown algebra explicitly.

\begin{eg} (\cite[2.3, 2.4]{G:struct}, cf.~\cite[1.9]{A:skew}) \label{browneg}
The principal examples of Brown $F$-algebras are denoted by $\B(J,
\D)$ for $J$ an Albert $F$-algebra and $\D$ a quadratic \'etale
$F$-algebra.  We set $\B(J, F \times F)$ to be the $F$-vector space
$\stbtmat{F}{J}{J}{F}$ with multiplication given by
\[
\tbtmat{\alpha_1}{j_1}{j'_1}{\beta_1}
\tbtmat{\alpha_2}{j_2}{j'_2}{\beta_2} = 
\tbtmat{\alpha_1 \alpha_2 + T(j_1, j'_2)}{\alpha_1 j_2 + \beta_2 j_1 +
j'_1 \times j'_2}{\alpha_2 j'_1 + \beta_1 j'_2 + j_1 \times
j_2}{\beta_1 \beta_2 + T(j_2, j'_1)},
\]
where $T$ is the bilinear trace form on $J$ and $\times$ is the
Freudenthal cross product (see \cite[\S 38]{KMRT} for more information about
these maps). 
The involution $-$ on $\B(J, F \times F)$ is given by
\[
\overline{\tbtmat{\alpha}{j}{j'}{\beta}} =
\tbtmat{\beta}{j}{j'}{\alpha}.
\]

The map $\varpi$ defined by
\[
\varpi \basmat = \tbtmat{\beta}{j'}{j}{\alpha}
\]
is an automorphism of $\B(J, F \times F)$ as an algebra with
involution.  For $\Delta$ a quadratic field extension of $F$ and
$\iota$ the nontrivial $F$-automorphism of $\D$, we define $\B(J, \D)$ 
to be the $F$-subspace of $\B(J, F \times F) \otimes_F \D$ fixed by
the $\iota$-semilinear automorphism $\varpi \otimes \iota$.  Note that
this construction is compatible with scalar extension in that for any
field extension $K$ of $F$, $\B(J, \D) \otimes_F K \cong \B(J
\otimes_F K, \D \otimes_F K)$.
\end{eg}

\begin{const} \label{quatconst}
Suppose that $K$ is a quadratic field extension of $F$, $Q$ is a
quaternion algebra over $F$ which is split by $K$, and $J$ is
an Albert $F$-algebra which is also split by $K$.
Then there is a gift $\G := (M_{28}(Q), \s, \pi)$
such that
\begin{enumerate}
\item $\Iso(\G)$ is split by $K$.
\item $\Der(\G)$ is isomorphic to the Lie algebra resulting from the
  Tits construction with $Q$ and $J$ as inputs.
\item If $Q$ is split, then $\G \cong \End{}(\B(J,K))$.
\item For $\gamma$ the unique symplectic involution on $Q$, the
  involution $\s$ is adjoint to the $\gamma$-hermitian form
  $\qform{1} \perp T$ for $T$ the trace on
  $J$.
\end{enumerate}
\end{const}

A thorough description of the Tits construction can be found in
\cite[\S 10]{Jac:ex}.

\newcommand{\psihat}{\underline{\psi}}
\newcommand{\phihat}{\underline{\varphi}}

\medskip
\begin{pf}
Consider the Brown algebras $\B := \B(J, K)$ and $\Bq := \B(J^d, K)$.
Their spaces of skew-symmetric elements are both spanned by some
element $s_0$ such that $F[s_0] \cong K$, and they are isomorphic as
algebras with involution over $K$, so $\B$ corresponds to a
1-cocycle $\psihat$ in $Z^1(K/F, \autp(\Bq))$, which must be
given by 
\[
\psihat_\iota \basmat = \tbtmat{\alpha}{\psi(j)}{\psi^\dagger(j')}{\beta}
\]
for some $\psi \in GL(J^d)(K)$ which preserves the norm on $J^d$.
Since $\psihat$ is a 
1-cocycle, $\psihat_\iota\ ^\iota \psihat_\iota = \Id_\B$, so we note
that $\psi \iota \psi^\dagger \iota = \Id_{J^d}$ for $\psi^\dagger \in
GL(J^d)(K)$ the unique map such that $T(\psi(j), \psi^\dagger(j')) =
T(j, j')$ for all $j, j' \in J^d$.

Since $Q$ is split by $K$, it is isomorphic to a quaternion algebra
$(a,b)_F$ for some $a, b \in \Fx$ such that $K = F(\sqrt{a})$.
Consider $t \in \EndF(\Bq)$ given by
\[
t \basmat = \tbtmat{\alpha /b}{b \psi
(j)}{\psi^\dagger (j')}{b^2 \beta}.
\]
By \cite[5.10]{G:struct}, the split triple system $\Md$ is the unique
triple system associated  
to $\Bq$.  The map $t$ is a similarity of $\Md$ with multiplier $b$, so
$\Int(t)$ is an automorphism of the gift associated to
$\End{}(\Bq)$, which is the split gift $\G^d$.  
Also,
\[
t \iota t \iota (x) = b x
\]
for all $x \in \Bq$.  Thus setting
$\eta_\iota = \Int(t)$ defines a 1-cocycle $\eta \in Z^1(K/F, \aut(\G^d))$,
which defines a nondegenerate gift $\bgift$ over $F$ which is split over $K$.

The map $H^1(K/F, \aut(\G^d)) \longto H^2(K/F, \mmu2)$
sends the class of $\G$ to the class of $A$ (which is the Tits algebra 
associated to $\Iso(\G)$), where we are identifying $H^2(K/F, \mmu2)$
and the subgroup of the Brauer group of $F$ consisting of algebras
which are split by $K$.
The image of $\eta$ under this
map is the 2-cocycle $f$ given by $f_{\iota, \iota} = b$.  This 2-cocycle
determines the quaternion algebra $Q$ by \cite[pp.~250, 251]{Sp:eq},
so $\eta$ determines a gift $\G$ of the form $(M_{28}(Q), \s, \pi)$.

The Lie algebras resulting from the Tits construction where $Q$ and
$J$ are split by $K$ are described in terms of Galois descent in
\cite[p.~87]{Jac:ex}, which immediately gives (2).  If $Q$ is split,
we may set $b = 1$, so part (3) is clear.  Thus we are left with
proving part (4).

For the Freudenthal triple system $(V, s, t)$ associated to $\B(J^d,
K)$ over $K$ we take $V = \stbtmat{F}{J^d}{J^d}{F} \otimes_F K$ and
\[
s\left( \stbtmat{\alpha_1}{j_1}{j'_1}{\beta_1},
\stbtmat{\alpha_2}{j_2}{j'_2}{\beta_2} \right) = \sqrt{a} \left
[ (\alpha_1 \beta_2 - \alpha_2 \beta_1) + (T(j_1, j'_2) - T(j'_1,
j_2)) \right].
\]
Note that $s$ is defined over $F$, where the action of $\iota$ on $V$
is by $\varpi \otimes \iota$.  This formula for the skew-symmetric
bilinear form comes from \cite[pp.~192, 195]{AF:CD}, where we have taken 
\[
s_0 := \stbtmat{\sqrt{a}}{}{}{-\sqrt{a}} \in \B(J^d, K) \otimes_F K.
\]
One can also find an explicit formula for the quartic form $q$ there,
which then implicitly specifies $t$, but we won't be using that.

Consider the embedding $\aut(\End{}(\B(J^d, K))) \injects \PSp(V,
s)$.  Since all nondegenerate skew-symmetric forms over $V$ of the same
dimension are isomorphic, there is some $\phihat \in \aut(V,s)(K)$ such
that
\[
\phihat\basmat =
\tbtmat{\alpha}{\varphi(j)}{\varphi^\dagger(j')}{\beta}
\]
and $\phihat^{-1} \psihat_\iota \,{^\iota \phihat} = \Id_V$.  Then
\[
\phihat^{-1} \eta_\iota \, {^\iota \phihat} \basmat = \tbtmat{\alpha /
b}{bj}{j'}{b^2 \beta}.
\]
Since $s$ is the skew-symemtric bilinear form constructed from
$\qform{1}$ with $c = b^{-1}$ and $T$ with $c = b$, by Lemma
\ref{symplem} we get that $\s$ is adjoint to $\qform{b^{-1}} 
\perp \qform{b}T \cong \qform{b} (\qform{1} \perp T)$, which is
similar to $\qform{1} \perp T$.
\end{pf}

\begin{eg} \label{isoeg}
Let $Q$ be a nonsplit quaternion algebra over a field $F$ with
splitting field $K$ a quadratic extension of $F$, and set $\G$ to be
the gift over $F$ constructed in \ref{quatconst} from $K$, $Q$, and
the split Albert $F$-algebra $J^d$.  Then $\G$ contains a rank 6
maximal singular ideal corresponding to the 6-dimensional maximal
singular ideal in $\B(J^d, \D)$ given in
\cite[7.6]{G:struct},
so by the description of the flag varieties from Section \ref{flagsec} the
non-end vertex of the length 2 arm of the Dynkin diagram of $\Iso(\G)$ 
is circled.  Since the Tits algebra of $\Iso(\G)$ is Brauer-equivalent 
to $Q$ and hence nonsplit, the end vertex of the long arm is not
circled, so $\Iso(\G)$ must be of type $E^9_{7,4}$ in the notation of
\cite[p.~59]{Ti:Cl}. 
\end{eg}

\section{Groups of type $E_7$ over real-closed fields} \label{galsec}

It is well-known that the Tits construction produces all Lie algebras
of type $E_7$ over a real-closed field $R$ (i.e., $R$ is real and
$R(\sqrt{-1})$ is algebraically closed), see \cite[pp.~120,
121]{Jac:ex}.  There are four of them, and they are obtained by
letting the quaternion algebra in the construction be split ($M_2(R)$)
or nonsplit (where $Q = H := (-1, -1)_R$, the unique nonsplit
quaternion algebra) and the Albert algebra $J$ be split ($J^d$) or
``minimally split'', so that $J = \h_3(O, 1)$, where $O$ is the Cayley
division algebra over $R$.  (For a definition and a general discussion of Cayley algebras,
please see \cite[Ch.~III, \S 4]{Schfr} or \cite[\S 33.C]{KMRT}.)
These quaternion and Albert algebras also
fill the hypotheses of our construction, and so our construction
produces all the groups of type $E_7$ over $R$.

For each of these four groups, we write $G = \Iso(M_{28}(Q), \s, \pi)$
and would like to say something about $\s$.  Since $\s$ is
symplectic, it is hyperbolic whenever $Q$ is split, so there is only
an issue when $Q = H$.  We have a very explicit
description of $\s$ coming from the construction, and here we use that
to calculate a particular invariant of $\s$, the {\em Witt index}
$w(M_{28}(Q), \s)$.
We know that $\s$ is adjoint to the hermitian form $h = \qform{1} \perp T$
over some $28$-dimensional $H$-vector space $V$ where $T$ is the trace
form of the Albert algebra $J$ in the construction, and the Witt index
is defined to be the number $2n$ such that $h$ is isomorphic to an
anisotropic form plus $n$ hyperbolic planes.  Equivalently, it is the
maximum rank of an isotropic right ideal $I$ 
of $(M_{28}(H), \s)$. To compute the Witt index, we apply a trick to
reduce to quadratic forms: We consider $V$ as a $(4 \cdot
28)$-dimensional vector space over $R$ endowed with a quadratic form
$q$ defined by $q(v) := h(v,v)$.  We say that $q$ is the {\em
quadratic trace form} of $h$.  This map (for general $V$) is an
injection from the Witt group of hermitian forms over $H$ to the Witt
group of quadratic forms over $R$ \cite[10.1.7]{Sch}, and in particular
a single hyperbolic plane has trace form four hyperbolic planes and
$\qform{1}$ has trace form
the reduced norm of $H$, which is $4\qform{1}$.  So the Witt
index of $\s$ is just half of the 
usual Witt index of the quadratic trace form $q$.  Since $q \cong 4
(\qform{1} \perp T)$ and our base field is real-closed, the Witt index
of $\s$ is twice that of the Witt index of the the quadratic form
$\qform{1} \perp T$.

Looking up the description of the trace form on $J$ from
\cite[37.9]{KMRT}, we find that the trace form on $J^d$ is 
$3\qform{1} \perp 12\hyp$ for $\hyp$ denoting a hyperbolic plane, and
the trace form on $\h_3(O, 1)$ is $27\qform{1}$. 
Thus we have the 
following table:
\begin{neqn}
\begin{tabular}{cc|ccc}
$Q$ & $J$ & $w(M_{28}(Q), \s)$ & type & comments \\ \hline
$M_2(R)$ & $J^d$ & 28 & $E^0_{7,7}$ & split \\
$M_2(R)$ & $\h_3(O, 1)$ & 28 & $E^{28}_{7,3}$ & \\
$H$ & $J^d$ & 24 & $E^9_{7,4}$ & \\
$H$ & $\h_3(O, 1)$ & 0 & $E^{133}_{7,0}$ & anisotropic/compact 
\end{tabular}
\end{neqn}
The entries in the ``type'' column refer to the types from Tits'
classification of simple groups \cite[p.~59]{Ti:Cl}.  The association
of types with choices of $Q$ and $J$ follow from the classification in
\cite{Jac:ex} and Example \ref{isoeg} or the fact that one can read off
from the Tits indices  
in \cite{Ti:Cl} that the groups of type $E^{28}_{7,3}$ have trivial Tits
algebra.

\begin{borel*}
One immediate consequence of this has to do with the Galois cohomology
of these groups.  Specifically, let $G = \Iso(M_{28}(H), \s, \pi)$ be
a simple simply connected 
group of type $E_7$ over $R$ and let $\Gbar = \aut(M_{28}(H), \s,
\pi)$ be its associated adjoint group.  Clearly $G$ is a subgroup of
$\Sp(M_{28}(H), \s)$ and $\Gbar$ is a subgroup of $\PSp(M_{28}(H),
\s)$, and the short exact sequence $1 \longto Z(G) \longto G \longto
\Gbar \longto 1$ and the corresponding one for the symplectic groups
induce a commutative diagram
\[
\begin{CD}
\Gbar(R) @>{\partial}>> H^1(R, \mmu2) @= \Rsq \\
@VVV @| @|\\
\PSp(M_{28}(H), \s) @>>> H^1(R, \mmu2) @= \Rsq.
\end{CD}
\]
That is, the image of $\partial$ is contained in the image of the
bottom arrow, which is known \cite[p.~425]{KMRT}.  It is the group of
similarity factors of the hermitian form $h$ (the hermitian form for
which $\s$ is adjoint), i.e., the elements $\mu \in \Fx$ such that
$\mu h \cong h$.
The construction of the quadratic trace form $q$
from $h$ makes it clear that any similarity factor of $h$ must also be
a similarity factor of $q$, so since $q$ is not hyperbolic and our
ground field is real-closed, the image of $\partial$ must be trivial.
\end{borel*}

\begin{borel*}
In the remaining two cases (where $Q$ is split), $\partial$ is
certainly surjective.  In fact, it is surjective over an arbitrary
field $F$ whenever $G$ is isomorphic to $\Inv(\M(J))$ for some Albert
$F$-algebra $J$.  This is because $\M(J)$ has a similarity with
multiplier $\lambda$ for all $\lambda \in \Fx$ by the construction in
\cite[p.~327]{Ferr:strict}.
\end{borel*}

\section*{Acknowledgements}
I would like to thank a referee for their helpful suggestions.

\providecommand{\bysame}{\leavevmode\hbox to3em{\hrulefill}\thinspace}

\bigskip

%
%
\noindent R.~Skip Garibaldi\\%
e-mail: {\tt skip@member.ams.org} \\%
web: {\tt http://www.math.ucla.edu/\~{}skip/}\\%

\noindent UC Los Angeles\\%
Dept.~of Mathematics\\%
Los Angeles, CA 90095-1555

\end{document}